\documentclass[aos,MSNbibl,dvips]{arximspdf}
\usepackage{graphicx}

\doi{10.1214/12-AOS1001} 
\referstodoi{10.1214/11-AOS949}
\volume{40}
\issue{4}
\pubyear{2012}
\firstpage{1997}
\lastpage{2004}

\makeatletter
\newcommand{\tr}{\operatorname{trace}}
\renewcommand{\epsilon}{\varepsilon}
\makeatother

\begin{document}
\begin{frontmatter}

\title{Discussion: Latent variable graphical model selection via
convex optimization}
\runtitle{Comment}

\begin{aug}
\author[A]{\fnms{Emmanuel J.} \snm{Cand\'es}\corref{}\ead[label=e1]{candes@stanford.edu}}
\and
\author[A]{\fnms{Mahdi} \snm{Soltanolkotabi}\ead[label=e2]{mahdisol@stanford.edu}}
\runauthor{E. J. Cand\'es and M. Soltanolkotabi}
\affiliation{Stanford University}
\address[A]{Department of Statisitcs \\
390 Serra Mall\\
Stanford University\\
Stanford, California 94305\\
USA\\
\printead{e1}\\
\phantom{E-mail:\ }\printead*{e2}} 
\end{aug}

\received{\smonth{4} \syear{2012}}



\end{frontmatter}

We wish to congratulate the authors for their innovative contribution,
which is bound to inspire much further research. We find latent
variable model selection to be a fantastic application of matrix
decomposition methods, namely, the superposition of low-rank and
sparse elements. Clearly, the methodology introduced in this paper is
of potential interest across many disciplines. In the following, we
will first discuss this paper in more detail and then reflect on the
versatility of the low-rank${}+{}$sparse decomposition.

\section*{Latent variable model selection}

The proposed scheme is an extension of the \emph{graphical lasso} of
Yuan and Lin \cite{Yuan07} (see also \cite{Banerjee,Friedman}), which
is a popular approach for learning the structure in an undirected
Gaussian graphical model. In this setup, we assume we have independent
samples $X \sim\mathcal{N}(0,\Sigma)$ with a covariance matrix $\Sigma$
exhibiting a sparse dependence structure but otherwise
unknown; that is to say, most pairs of variables are conditionally
independent given all the others. Formally, the concentration matrix
$\Sigma^{-1}$ is assumed to be sparse. A natural fitting procedure is
then to regularize the likelihood by adding a term proportional to the
$\ell_1$ norm of the estimated inverse covariance matrix $S$:
%
\begin{equation}\label{eq:glasso}
\mbox{minimize}-\!\ell(S,\Sigma_0^n) + \lambda\|S\|_1
\end{equation}
under the constraint $\mathbf{S}\succeq 0$, where $\Sigma_0^n$ is the
empirical covariance matrix and $\|S\|_1 = \sum_{ij}
|S_{ij}|$. (Variants are possible depending upon whether or not one
would want to penalize the diagonal elements.) This problem is convex.

When some variables are unobserved---the observed and hidden variables
are still jointly Gaussian---the model above may not be appropriate
because the hidden variables can have a confounding effect. An example
is this: we observe stock prices of companies and would like to infer
conditional (in)dependence. Suppose, however, that all these companies
rely on a commodity, a source of energy, for instance, which is not
observed. Then the stock prices might appear dependent even though
they may not be once we condition on the price of this commodity. In
fact, the marginal inverse covariance of the observed variables
decomposes into two terms. The first is the concentration matrix of
the observed variables in the full model conditioned on the latent
variables. The second term is the effect of marginalization over the
hidden variables. Assuming a sparse graphical model, the first term is
sparse, whereas the second term may have low rank; in particular, the
rank is at most the number of hidden variables. The authors then
penalize the negative log-likelihood with a term proportional to
%
\begin{equation}\label{eq:L+S}
\gamma\|S\|_1 + \tr(L)
\end{equation}
since the trace functional is the usual convex surrogate for the rank
over the cone of positive semidefinite matrices. The constraints are
$S \succ L \succeq0$.

\section*{Adaptivity} The penalty (\ref{eq:L+S}) is simple and flexible
since it does not really make special parametric assumptions. To be
truly appealing, it would also need to be adaptive in the following
sense: suppose there is no hidden variable, then does the low-rank${}+{}$sparse model
(L${}+{}$S) behave as well or nearly as well as the graphical
lasso? When there are few hidden variables, does it behave nearly as
well? Are there such theoretical guarantees? If this is the case, it
would say that using the L${}+{}$S model would protect against the danger of
not having accounted for all possible covariates. At the same time, if
there were no hidden variable, one would not suffer any loss of
performance. Thus, we would get the best of both worlds.

At first sight, the analysis presented in this paper does not allow us
to reach this conclusion. If $X$ is $p$-dimensional, the number of
samples needed to show that one can obtain accurate estimates scales
like $\Omega(p/\xi^4)$, where $\xi$ is a modulus of continuity
introduced in the paper that is typically much smaller than 1. We can
think of $1/\xi$ as being related to the maximum degree $d$ of the
graph so that the condition may be interpreted as having a number of
observations very roughly scaling like $d^4 p$. In addition, accurate
estimation holds with the proviso that the signal is strong enough;
here, both the minimum nonzero singular value of the low-rank
component and the minimum nonzero entry of the sparse component scale
like $\Omega(\sqrt{p/n})$. On the other hand, when there are no hidden
variables, a line of work \cite{Meinshausen06,Ravikumar11,Rothman08}
has established that we could estimate the concentration matrix with
essentially the same accuracy if $n = \Omega(d^2 \log p)$
and the magnitude of the minimum nonvanishing value of the
concentration matrix scales like $\Omega(\sqrt{n^{-1} \log p})$. As
before, $d$ is the maximum degree of the graphical model. In the
high-dimensional regime, the results offered by this literature seem
considerably better. It would be interesting to know whether this
could be bridged, and if so, under what types of conditions---if
any.\vadjust{\goodbreak}

Interestingly, such adaptivity properties have been established for
related problems. For instance, the L${}+{}$S model has been used to suggest
the possibility of a principled approach to robust principal component
analysis \cite{Candes11}. Suppose we have incomplete and corrupted
information about an $n_1 \times n_2$ low-rank matrix $L^0$. More
precisely, we observe $M_{ij} = L^0_{ij} + S^0_{ij}$, where $(i,j) \in
\Omega_{\mathrm{obs}} \subset\{1, \ldots, n_1\} \times\{1, \ldots,
n_2\}$. We think of $S^0$ as a corruption pattern so that some entries
are totally unreliable but we do not know which ones. Then
\cite{Candes11} shows that under rather broad conditions, the solution
to
\begin{eqnarray}\label{eq:RPCA}
&&\mbox{minimize }  \|L\|_* + \lambda\|S\|_1\nonumber\\ [-8pt]\\ [-8pt]
&&\quad\mbox{subject to }   M_{ij} = L_{ij} + S_{ij},  (i,j) \in
\Omega_{\mathrm{obs}},\nonumber
\end{eqnarray}
where $\|L\|_*$ is the nuclear norm, recovers $L^0$ exactly. Now
suppose there are no corruptions. Then we are facing a matrix
completion problem and, instead, one would want to minimize the
nuclear norm of $L$ under data constraints. In other words, there is
no need for $S$ in (\ref{eq:RPCA}). The point is that there is a
fairly precise understanding of the minimal number of samples needed
for this strategy to work; for incoherent matrices \cite{CR08},
$|\Omega_{\mathrm{obs}}|$ must scale like $(n_1 \vee n_2) r \log^2 n$,
where $r$ is the rank of $L^0$. Now some recent work~\cite{Li}
establishes the adaptivity in question. In details, (\ref{eq:RPCA})
recovers $L^0$ from a minimal number of samples, in the sense defined
above, even though a positive fraction may be corrupted. That is, the
number of reliable samples one needs, regardless of whether corruption
occurs, is essentially the same. Results of this kind extend to other
settings as well. For instance, in sparse regression or compressive
sensing we seek a sparse solution to $y = X b$ by minimizing the
$\ell_1$ norm of $b$. Again, we may be worried that some equations are
unreliable because of gross errors and would solve, instead,
\begin{eqnarray}\label{eq:RCS}
&&\mbox{minimize }   \|b\|_1 + \lambda\|e\|_1\nonumber\\ [-8pt]\\ [-8pt]
&&\quad\mbox{subject to }    y = Xb + e\nonumber
\end{eqnarray}
to achieve robustness. Here, \cite{Li} shows that the minimal number
of reliable samples/equations required, regardless of whether the data
is clean or corrupted, is essentially the same.



\section*{The versatility of the L${}+{}$S model}

We now move to discuss the L${}+{}$S model more generally and survey a set
of circumstances where it has proven useful and powerful.
To begin with, methods which simply minimize an $\ell_1$ norm, or a
nuclear norm, or a combination thereof are seductive because they are
flexible and apply to a rich class of problems. The L${}+{}$S model is
nonparametric and does not make many assumptions. As a result, it is
widely applicable to problems ranging from latent variable model
selection \cite{Chandrasekaran11} (arguably one of the most subtle and beautiful applications
of this method) to video surveillance in computer vision and document
classification in machine learning \cite{Candes11}. In any given application, when
much is known about the problem, it may not return the best possible
answer, but our experience is that it is always fairly competitive.
That is, the little performance loss we might encounter is more than
accounted for by the robustness we gain vis a vis various modeling
assumptions, which may or may not hold in real applications. A few
recent applications of the L${}+{}$S model demonstrate its flexibility and
robustness.



\section*{Applications in computer vision}
The L${}+{}$S model has been applied
to address several problems in computer vision, most notably by the
group of Yi Ma and colleagues. Although the low-rank${}+{}$sparse model may
not hold precisely, the nuclear${}+{}$$\ell_1$ relaxation appears
practically robust. This may be in contrast with algorithms which use
detailed modeling assumptions and may not perform well under slight
model mismatch or variation.

\subsection*{Video surveillance} An important task in computer vision is to
separate background from foreground. Suppose we stack a sequence of
video frames as columns of a matrix (rows are pixels and columns time
points), then it is not hard to imagine that the background will have
low-rank since it is not changing very much over time, while the
foreground objects, such as cars, pedestrians and so on, can be seen as
a sparse disturbance. Hence, finding an L${}+{}$S decomposition offers a new
way of modeling the background (and foreground). This method has been
applied with some success \cite{Candes11}; see also the online videos
\href{http://www.youtube.com/watch?feature=player_embedded&v=RPmr8WLkBSo}{Video
1} and
\href{http://www.youtube.com/watch?feature=player_embedded&v=Yxj1_52EAXA}{Video
2}.

\subsection*{From textures to 3D} 
\begin{figure}
\begin{tabular}{@{}cc@{}}

\includegraphics{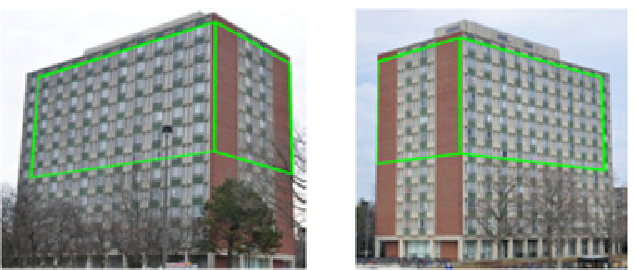}
 & \includegraphics{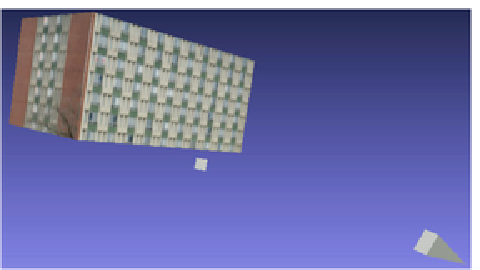}\\
(a) & (b)
\end{tabular}
\caption{(\textup{a}) Pair of images from distinct viewpoints. (\textup{b}) 3D
reconstruction (TILT) from photographs in (\textup{a}) using the L${}+{}$S
model. The geometry is recovered from two images.}\label{fig:tilt}
\end{figure}

\begin{figure}[b]

\includegraphics{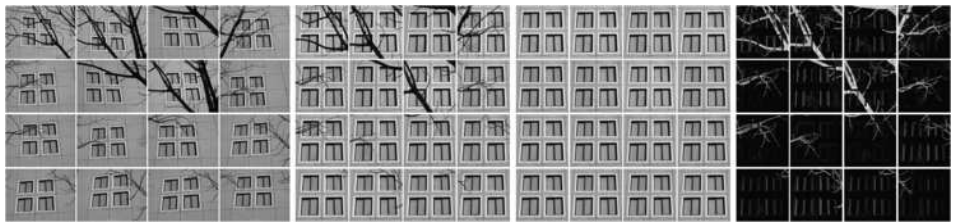}

\caption{We are given the 16 images on the right. The task is to
remove the clutter and align the images. Stacking each image as a
column of a matrix, we look for planar homeographies that reveal a
low-rank plus sparse structure \protect\cite{Peng10}. From left to right:
original data set, aligned images, low-rank component (columns of
$L$), sparse component (columns of $S$).}\label{fig:rasl}
\end{figure}

One of the most fundamental
steps in computer vision consists of extracting relevant features that
are subsequently used for high-level vision applications such as 3D
reconstruction, object recognition and scene understanding. There has
been limited success in extracting stable features across variations
in lightening, rotations and viewpoints. Partial occlusions further
complicate matters. For certain classes of 3D objects such as images
with regular symmetric patterns/textures, one can bypass the
extraction of local features to recover 3D structure from 2D views.
To fix ideas, a vertical or horizontal strip can be regarded as a
rank-1 texture and a corner as a rank-2 texture. Generally speaking,
surfaces may exhibit a low-rank texture when seen from a suitable
viewpoint; see Figure~\ref{fig:tilt}. However, their 2D projections as
captured by a camera will typically not be low rank. To see why,
imagine there is a low-rank texture ${L}^0(x,y)$ on a planar
surface. The image we observe is a transformed version of this
texture, namely, $L^0\circ\tau^{-1}(x,y)$. A technique named TILT
\cite{Zhang12} recovers $\tau$ simply by seeking a low-rank and sparse
superposition. In spite of idealized assumptions, Figures~\ref{fig:tilt} and \ref{fig:rasl} show that the L${}+{}$S model works well
in practice.

%
%
%
%

%

\subsection*{Compressive acquisition} In the spirit of compressive sensing,
the L${}+{}$S model can also be used to speed up the acquisition of large
data sets or lower the sampling rate. At the moment, the theory of
compressive sensing relies on the sparsity of the object we wish to
acquire, however, in some setups the L${}+{}$S model may be more
appropriate. To explain our ideas, it might be best to start with two
concrete examples. Suppose we are interested in the efficient
acquisition of either (1) a hyper-spectral image or (2) a video
sequence. In both cases, the object of interest is a data matrix $M$
which is $N \times d$, where each column is an $N$-pixel image and each
of the $d$ columns corresponds to a specific wavelength (as in the
hyper-spectral example) or frame (or time point as in the video
example). In the first case, the data matrix may be thought of as
$M(x,\lambda)$, where $x$ indexes position and $\lambda$ wavelength,
whereas in the second example, we have $M(x,t)$ where $t$ is a time
index. We would like to obtain a sequence of highly resolved images
from just a few measurements; an important application concerns
dynamic magnetic resonance imaging where it is only possible to
acquire a few samples in $k$-space per time interval.

Clearly, frames in a video sequence are highly correlated in time. And
in just the same way, two images of the same scene at nearby
wavelengths are also highly correlated. Obviously, images are
correlated in space as well. Suppose that $W \otimes F$ is a tensor
basis, where $W$ sparsifies images and $F$ time traces ($W$ might be a
wavelet transform and $F$ a Fourier transform). Then we would expect
$WMF$ to be a nearly sparse matrix. With undersampled data of the form
$ y = \mathcal{A}(M) + z$, where $\mathcal{A}$ is the operator supplying
information about $M$ and $z$ is a noise term, this leads to the
low-rank${}+{}$sparse decomposition problem
\begin{eqnarray}\label{eq:prop1}
&&\mbox{minimize }   \|X\|_* + \lambda\|W X F\|_1\nonumber\\ [-8pt]\\ [-8pt]
&&\quad\mbox{subject to }   \|\mathcal{A}(X) - y\|_2 \le\epsilon,\nonumber
\end{eqnarray}
where $\epsilon^2$ is the noise power. A variation, which is more in
line with the discussion paper is a model in which $L$ is a low-rank
matrix modeling the static background, and $S$ is a sparse matrix
roughly modeling the innovation from one frame to the next; for
instance, $S$ might encode the moving objects in the foreground. This
would give
\begin{eqnarray}\label{eq:prop2}
&&\mbox{minimize }  \lambda\|L\|_* + \|W S F\|_1\nonumber\\ [-8pt]\\ [-8pt]
&&\quad\mbox{subject to }  \|\mathcal{A}(L+S) - y\|_2 \le\epsilon.\nonumber
\end{eqnarray}
One could imagine that these models might be useful in alleviating the
tremendous burden on system resources in the acquisition of ever
larger 3D, 4D and 5D data sets.

We note that proposals of this kind have begun to emerge. As we were
preparing this commentary, we became aware of \cite{Golbabaee12},
which suggests a model similar to~(\ref{eq:prop1}) for hyperspectral
imaging. The difference is that the second term in (\ref{eq:prop1}) is
of the form $\sum_i \|X_i\|_{\mathrm{TV}}$ in which $X_i$ is the $i$th
column of~$X$, the image at wavelength $\lambda_i$; that is,~we minimize
the total variation of each image, instead of looking for sparsity
simultaneously in space and wavelength/frequency. The results in
\cite{Golbabaee12} show that dramatic undersampling ratios are
possible. In medical imaging, movement due to respiration can degrade
the image quality of Computed Tomography (CT), which can lead to
incorrect dosage in radiation therapy. Using time-stamped data, 4D CT
has more potential for precise imaging. Here, one can think of the
object as a matrix with rows labeling spatial variables and columns
time. In this context, we have a low-rank (static) background and a
sparse disturbance corresponding to the dynamics, for example,~of the
heart in
cardiac imaging. The recent work \cite{Gao11} shows how one can use
the L${}+{}$S model in a fashion similar to~(\ref{eq:prop2}). This has
interesting potential for dose reduction since the approach also
supports substantial undersampling.

\section*{Connections with theoretical computer science and future
directions} A~class of problems where further study is required
concerns situations in which the low-rank and sparse components have a
particular structure. One such problem is the \emph{planted clique
problem}. It is well known that finding the largest clique in a
graph is NP hard; in fact, it is even NP-hard to approximate the size
of the largest clique in an $n$ vertex graph to within a factor
$n^{1-\epsilon}$. Therefore, much research has focused on an
``easier'' problem. Consider a random graph $G(n,1/2)$ on $n$ vertices
where each edge is selected independently with probability $1/2$. The
expected size of its largest clique is known to be $(2-o(1))\log n$.
The planted clique problem adds a clique of size $k$ to $G$. One hopes
that it is possible to find the planted clique in polynomial time
whenever $k \gg\log n$. At this time, this is only known to be
possible if $k$ is on the order of $\sqrt{n}$ or larger. In spite of
its seemingly simple formulation, this problem has eluded theoretical
computer scientists since 1998, and is regarded as a notoriously
difficult problem in modern combinatorics. It is also fundamental to
many areas in machine learning and pattern recognition. To emphasize
its wide applicability, we mention a new connection with game
theory. Roughly speaking, the recent work \cite{Hazan10} shows that
finding a near-optimal Nash equilibrium in two-player games is as hard
as finding hidden cliques of size $k = C_0   \log n$, where $C_0$ is
some universal
constant. 

One can think about the planted clique as a low rank${}+{}$sparse
decomposition problem. To be sure, the adjacency matrix of the graph
can be written as the sum of two matrices: the low-rank component is
of rank 1 and represents the clique of size $k$ (a submatrix with all
entries equal to $1$); the sparse component stands for the random
edges (and with $-1$ on the diagonal if and only if that vertex
belongs to the hidden clique). Interestingly, low-rank${}+{}$sparse
regularization based on nuclear and $\ell_1$ norms have been applied
to this problem \cite{Doan10}. (Here the clique is both low-rank and
sparse and is the object of interest so that we minimize $\|X\|_* +
\lambda\|X\|_1$ subject to data constraints.)
These proofs show that these methods find cliques of size
$\Omega(\sqrt{n})$, thus recovering the best known results, but they
may not be able to break this barrier. It is interesting to
investigate whether tighter relaxations, taking into account the
specific structure of the low-rank and sparse components, can do
better.


%

\printaddresses

\end{document}